\documentclass[preprint,12pt]{elsarticle}

\usepackage{amsmath}
\usepackage{amssymb}
\usepackage{amscd}
\usepackage{epsfig}
\usepackage{color}
\usepackage{subfigure}
\usepackage{graphicx}

\newcommand{\N}{\mathbb N}
\newcommand{\bee}{\begin{equation}}
\newcommand{\eee}{\end{equation}}

\newcommand{\Lb}{\mbox {\boldmath ${\Lambda}$}}

\newcommand{\Lbs}{\mbox{\scriptsize\boldmath ${\Lambda}$}}

\def\ov{\overline}

\def\b0{{\bf 0}}

\newcommand{\Pb}{\mbox {\bf P}}

\definecolor{darkgreen}{rgb}{0.3,0.6,0.1}

\newcommand{\be}{\begin{eqnarray}}
\newcommand{\ee}{\end{eqnarray}}
\newcommand{\supp}{\mbox{\rm supp}}

\newcommand{\dens}{\mbox{\rm dens}}

\newcommand{\Vol}{\mbox{\rm Vol}}

\newcommand{\eps}{{\mbox{$\epsilon$}}}

\newcommand{\R}{{\mathbb R}}
\newcommand{\Q}{{\mathbb Q}}
\newcommand{\Z}{{\mathbb Z}}
\newcommand{\C}{{\mathbb C}}

\newcommand{\Ak}{{\mathcal A}}

\newcommand{\Dk}{{\mathcal D}}

\newcommand{\Sk}{{\mathcal S}}
\newcommand{\Tk}{{\mathcal T}}

\newcommand{\Lam}{{\Lambda}}

\newcommand{\lam}{\lambda}

\newtheorem{theorem}{Theorem}[section]
\newtheorem{lemma}[theorem]{Lemma}

\newtheorem{cor}[theorem]{Corollary}

\newtheorem{example}[theorem]{Example}
\newtheorem{remark}[theorem]{Remark}

\newtheorem{defi}[theorem]{Definition}
\numberwithin{equation}{section}

\journal{European Journal of Combinatorics}

\begin{document}

\begin{frontmatter}

\title{Overlap coincidence to strong coincidence in substitution tiling dynamics }

\author{Shigeki Akiyama $^{\,\rm a}$,  Jeong-Yup Lee $^{\,\rm b*}$ }

\address{
a: 
Institute of Mathematics, University of Tsukuba, 1-1-1 Tennodai, \\
\hspace*{2.5em} Tsukuba, Ibaraki, Japan (zip:305-8571);
{akiyama@math.tsukuba.ac.jp}

\smallskip
b: Dept. of Math. Edu., Kwandong University, 24, 579 Beon-gil, Beomil-ro, Gangneung, \\ 
\hspace*{2.5em} Gangwon-do, 210-701 Republic of Korea;
{jylee@kd.ac.kr, jeongyuplee@yahoo.co.kr}; \\
Tel. 82-33-649-7776

\smallskip
{*Corresponding author\\} }

\begin{abstract}
Overlap coincidence is an equivalent criterion to pure discrete spectrum of the dynamics of self-affine tilings in $\R^d$. 
In the case of $d=1$, strong coincidence on $m$-letter irreducible substitution has been introduced 
in \cite{Dek, AI} which implies that the system is 
metrically conjugate to a domain exchange in $\R^{m-1}$.
However being a domain exchange does not imply the property of pure discrete spectrum of the tiling dynamics. 
The relation between  two coincidences has not been established completely. 
In this paper we generalize  strong coincidence to higher dimensions
and show the implication from  overlap coincidence to the new strong coincidence
when the associated height group is trivial.
Furthermore we introduce a new criterion `simultaneous coincidence' and show the implication  from  overlap coincidence to the simultaneous coincidence. The triviality of  height group is shown in \cite{BK, SingThesis} for $1$-dimension irreducible Pisot substitutions. 
\end{abstract}


\end{frontmatter}

\section{Introduction}

The principal aim of this paper is to give a better understanding of pure discrete spectrum
of self-affine tiling dynamical systems which have zero-entropy and whose spectral type varies from weakly-mixing to 
pure discrete. The study is strongly motivated by  atomic configurations
of quasicrystals, which show pure point diffraction. Indeed,  equivalence of 
pure point diffraction of quasicrystal structure and  pure discrete spectrum of its 
associated dynamical system is 
known in quite a general setting. If we restrict ourselves to $1$-dimension 
substitutive systems, 
the problem of pure discrete spectrum can be reformulated using  notions of `coincidence' on word combinatorics.
Many notions of coincidences have been introduced in the study of the dynamical spectrum of self-affine tilings. 
A lot of coincidences among these are proved to be equivalent. 
However  relation between strong coincidence of $1$-dimension Pisot substitution 
and  other coincidences is not completely understood.

Overlap coincidence is defined in self-affine tiling and characterizes pure discrete spectrum of the tiling dynamics \cite{soltil, LMS2, Lee}.
Basically what it means is that every two tiles in the tiling, which overlap after translating by a return vector, have at least one tile in common after some iterations (see subsection \ref{Overlap-coincidence} for the detail). 
This coincidence has been proved to be equivalent with algebraic coincidence \cite{Lee} and super coincidence \cite{Siegel}. 
For an irreducible Pisot substitution in $1$-dimension, 
all of these coincidences are equivalent to the fact that 
balanced pair algorithm terminates and each balanced pair leads to a coincidence \cite{HS,SS}.
This is a combinatorial condition that we can quickly check for a given substitution.

Strong coincidence has been introduced by \cite{Dek} in constant length substitution sequences, generalized in unimodular Pisot substitutions by \cite{AI}, and extended in the case of non-unimodular Pisot substitutions by \cite{Siegel}. 
This combinatorial condition guarantees that there is a geometric realization of substitutions which is metrically conjugate to a domain exchange.
In  view of the balanced pair algorithm, 
 strong coincidence only implies that every balanced pair leads to a coincidence,
i.e., termination of the balanced pair algorithm seems to be necessary to assume.
In \cite{HS} Hollander and Solomyak proved the termination of the balanced pair algorithm in two letter case,
establishing the equivalence between strong coincidence and pure discrete spectrum.
From this equivalence with the result of Barge and Diamond \cite{BD} which guarantees the strong coincidence in two-letter case, we obtained that 
two-letter irreducible Pisot substitution dynamical systems have  pure discrete spectrum.

However apart from  two-letter irreducible 
Pisot substitution sequences, the relation between the strong coincidence and pure discrete spectrum is not clearly understood. 
This relation is important for an approach towards `Pisot substitution conjecture' \cite{ABBLS}.
Recently Nakaishi \cite{Nakaishi} claimed that the dynamics of irreducible unimodular Pisot substitution sequence satisfying strong coincidence has pure discrete spectrum through domain exchange flow. Though we do not yet have a full account of this claim, 
it indicates that strong coincidence in irreducible unimodular Pisot substitution sequences implies overlap coincidence in irreducible unimodular Pisot substitution tilings in $\R$.
This gives us a motivation to look at the other direction. 
As mentioned above, strong coincidence and overlap coincidence are 
both combinatorial objects which are defined in $1$-dimension substitution sequence. 
However, it seems that there is no reason to restrict it to $1$-dimension substitutions to see the relationship between 
these coincidences and seems to be better to transfer the problem into geometric setting.
So we first
generalize  strong coincidence in substitution sequences into substitution tiling in $\R^d$, and give a stronger version `simultaneous coincidence' in \S 2
and show that overlap coincidence implies simultaneous 
coincidence in Pisot family substitution tiling in $\R^d$
provided that the associated substitution Delone multi-color set is `admissible' (Def. \ref{admissible}) and the corresponding `height group' (Def. \ref{Height}) is trivial.
Since it is shown in \cite{BK, SingThesis} that this group is trivial for irreducible 
Pisot substitutions, we know that every irreducible Pisot substitution having 
pure discrete spectrum must admit the simultaneous coincidence. 

Here we point out that our higher dimensional 
generalization of strong coincidence is somewhat too weak condition by itself. It must be fulfilled in conjunction with some constraints to make a reasonable connection to other coincidences.
Indeed, without the constraints of  
admissibility and  trivial height group, for a given substitution tiling, 
it is always possible to build an associated substitution Delone 
multi-color set satisfying our extended notion of strong coincidence. 
The first author \cite{Ak} obtained a converse statement from strong coincidence to overlap coincidence using 
the extended notion of strong coincidence defined in this paper. 
But in \cite{Ak}, strong coincidence is required 
for many choices of admissible {\it control points} having trivial height group to deduce overlap coincidence.
It is of interest to minimize the constraints which are sufficient to derive 
 converse direction from strong coincidence to overlap coincidence.

\section{Preliminary}

\subsection{Tiles and tilings}

A {\em tile} in $\R^d$ is defined as a pair $T=(A,i)$ where $A=\supp(T)$
(the support of $T$) is a compact
set in $\R^d$ which is the closure of its interior, and
$i=l(T)\in \{1,\ldots,m\}$
is the color of $T$. We let $g+T = (g+A,i)$ for $g \in \R^d$. We say that
a set $P$ of tiles is a {\em patch} if the number of tiles in $P$ is
finite and the tiles of $P$ have mutually disjoint
interiors. A {\em tiling} of $\R^d$ is a set $\Tk$ of tiles such that $\R^d = \bigcup \{\supp(T) : T \in \Tk \}$ and distinct tiles have disjoint interiors.

\begin{defi}\label{def-subst}
{\em Let $\Ak = \{T_1,\ldots,T_m\}$ be a finite set of tiles in
$\R^d$ such that $T_i=(A_i,i)$; we will call them {\em
prototiles}. Denote by ${\mathcal{P}}_{\Ak}$ the set of non-empty
patches. 
We say that $\Omega: \Ak \to {\mathcal{P}}_{\Ak}$ is a {\em
tile-substitution} (or simply {\em substitution}) with a
$d\times d$ expansive matrix $Q$ if there exist finite sets $\Dk_{ij}\subset
\R^d$ for $i,j \le m$ such that
\begin{equation}
\Omega(T_j)=
\{u+T_i:\ u\in \Dk_{ij},\ i=1,\ldots,m\}
\label{subdiv}
\end{equation}
with
\begin{eqnarray} \label{tile-subdiv}
Q A_j = \bigcup_{i=1}^m (\Dk_{ij}+A_i) \ \ \ \mbox{for} \ j\le m.
\end{eqnarray}
Here all sets in the right-hand side must have disjoint interiors;
it is possible for some of the $\Dk_{ij}$ to be empty.}
\end{defi}

\noindent The substitution (\ref{subdiv}) is extended to all
translates of prototiles by $\Omega(x+T_j)= Q x + \Omega(T_j)$ and
to patches and tilings by $\Omega(P)=\bigcup\{\Omega(T):\ T\in
P\}$. The substitution $\Omega$ can be iterated, producing larger
and larger patches $\Omega^k(P)$.
We say that $\Tk$ is a {\em substitution tiling} if $\Tk$ is a
tiling and $\Omega(\Tk) = \Tk$ with some substitution $\Omega$. In
this case, we also say that $\Tk$ is a {\em fixed point} of
$\Omega$. We say that a substitution tiling is {\em primitive} if
the corresponding substitution matrix $S$, with $S_{ij}= \sharp
(\Dk_{ij})$, is primitive, and {\em irreducible} if the characteristic polynomial of $S$ is irreducible. 
We say that $\Tk$ has {\em finite local complexity} (FLC) 
if $\forall \ R > 0$, $\exists$ finitely many translational classes of patches whose support lies in some ball of radius $R$.
A tiling $\Tk$ is {\em repetitive} if for any compact set
$K \subset \R^d$, $\{t \in \R^d : \Tk \cap K = (t + \Tk) \cap K\}$ is relatively dense.
A repetitive fixed point of a primitive
tile-substitution with FLC is called a {\em self-affine tiling}.
Let $\lambda>1$ be the Perron-Frobenius eigenvalue  of the substitution matrix $S$.
Let
$ D =\{\lam_1,\ldots,\lam_{d}\}$ be the set of (real and complex)
eigenvalues of $Q$. We
say that $Q$ (or the substitution $\Omega$) fulfills a {\em Pisot family} if for every $\lam\in
D$ and every Galois conjugate $\lam'$ of $\lam$, if
$\lam'\not\in D$, then $|\lam'| < 1$.

\medskip

Given a tiling $\Tk$ in $\R^d$, we define the {\em tiling space} as the
orbit closure of $\Tk$ under the translation action: $X_{\Tk} =
\ov{\{-g + \Tk :\,g \in \R^d \}}$, in the well-known ``local
topology'': for a small $\eps>0$ two tilings $\Sk_1,\Sk_2$ are
$\eps$-close if $\Sk_1$ and $\Sk_2$ agree on the ball of radius
$\eps^{-1}$ around the origin, after a translation of size less
than $\eps$. Then $X_{\Tk}$ is compact and we get a topological dynamical system
$(X_{\Tk},\R^d)$ where $\R^d$ acts by translations. This system is
minimal (i.e.\ every orbit is dense) whenever $\Tk$ is repetitive.

\subsection{Point sets}

\noindent Recall that a Delone set is a
relatively dense and uniformly discrete subset of $\R^d$. We say
that $\Lb=(\Lambda_i)_{i\le m}$ is a {\em Delone multi-color set}
in $\R^d$ if each $\Lambda_i$ is Delone and
$\supp(\Lb):=\cup_{i=1}^m \Lambda_i \subset \R^d$ is Delone. 
A {\em cluster} of $\Lb$ is a family of points $\Pb = (P_i)_{i \le m}$ where $P_i \subset \Lam_i$ is finite for all $i \le m$.
We say that $\Lam \subset \R^d$ is a {\em Meyer set}
if it is a Delone set and $\Lam - \Lam$ is uniformly discrete in $\R^d$ \cite{Lag}.
The
colors of points in the Delone multi-color set have the same
meaning as the colors of tiles on tilings.
Various notions such as primitivity, FLC and repetitivity in point sets are defined in similar way as in tilings. 

\begin{defi} \label{def-subst-mul}
{\em $\Lb = (\Lam_i)_{i\le m}$ is called a {\em
substitution Delone multi-color set} if $\Lb$ is a Delone multi-color set and
there exist an expansive matrix
$Q$ and finite sets $\Dk_{ij}$ for $i,j\le m$ such that
\be \label{eq-sub}
\Lambda_i = \bigcup_{j=1}^m (Q \Lambda_j + \Dk_{ij}),\ \ \ i \le m,
\ee
where the unions on the right-hand side are disjoint.}
\end{defi}

For any given substitution Delone multi-colour set $\Lb = (\Lambda_i)_{i \le m}$,
we define $\Phi_{ij} = \{ f : x \mapsto Qx + a \, : \,a \in \Dk_{ij}\}$.
Then $\Phi_{ij}(\Lam_j) = Q \Lam_j + \Dk_{ij}$, where $i \le m$. We define $\Phi$ an $m \times m$ array for which each entry is $\Phi_{ij}$.
For any $k \in \Z_+$ and $x \in \Lam_j$ with $j \le m$, we let
$\Phi^k (x) =  \Phi^{k-1}((\Phi_{ij}(x))_{i \le m})$.
For any $k \in \Z_+$, $\Phi^k (\Lb) = \Lb$ and
$\Phi^k (\Lam_j) = \bigcup_{i \le m}(Q^k \Lam_j + (\Dk^k)_{ij})$ where
\[(\Dk^k)_{ij} = \bigcup_{n_1,n_2,\dots,n_{(k-1)} \le m}
(\Dk_{in_1} + Q \Dk_{n_1 n_2} + \cdots + Q^{k-1} \Dk_{n_{(k-1)} j}).
\]
We say that a cluster $\Pb$ is {\em legal} if it is a translate of a subcluster of
a cluster generated from one point of $\Lb$, that is to say, $a + \Pb \subset \Phi^k (x)$ (i.e. $a + \Pb \subset \bigcup_{i=1}^m (Q^k x + (\Dk^k)_{ij})$) for some $k \in \Z_+$, $a \in \R^d$, $x \in \Lam_j$, and $j \le m$.

\begin{defi}
{\em Let $\Lb$ be a primitive substitution Delone multi-color set in $\R^d$ 
with an expansive matrix $Q$. 
We say that $\Lb$ admits an {\em algebraic coincidence} if there exist $M \in \Z_+$ and $\xi \in \Lam_i$ for some $i \le m$ such that $Q^M \Xi(\Lb) \subset \Lam_i - \xi$, where $\Xi(\Lb) = \bigcup_{i \le m} (\Lam_i - \Lam_i)$}. 
\end{defi}

Since $\Lb$ is primitive, there exists $N \in \N$ such that $Q^N \Xi(\Lb) \subset \Lam_i - \Lam_i$ for any $i \le m$. So $\Lb$ admits an algebraic coincidence if and only if there exists $M \in \N$ and $\xi \in \Lam_i$ such that $Q^M (\Lam_i - \Lam_i) \subset \Lam_i - \xi$.
We say that $\Lb$ is a {\em Pisot family substitution Delone multi-color set} if $Q$ fulfills the condition of Pisot family.

\medskip

If a self-affine tiling $\Tk$ is
given, we can get an associated substitution Delone multi-color set
$\Lb_{\Tk} = (\Lam_i)_{i \le m}$ of $\Tk$ taking representative points of tiles in the relatively same positions for the same color tiles in the tiling (see \cite[Lemma\,5.4]{Lee}). 
On the other hand, if $\Lb$ is a primitive substitution Delone multi-color set for which every $\Lb$-cluster is legal, then $\Lb + \mathcal{A}:= \{x + T_i : x \in \Lam_i, i \le m \}$ is a tiling of $\R^d$, where $\mathcal{A}$ is the set of prototiles from the associated tile equations (see \cite{LMS2}).
This bijection establishes a topological
conjugacy of $(X_{\Lbs},\R^d)$ and $(X_{\Tk},\R^d)$. 

Note that if we translate representative points of tiles of $\Tk$ 
by $\Lam'_i=\Lam_i-g_i$, then the set equation will be
$$
\Lam'_i= \bigcup_{j=1}^m Q \Lam'_j + \mathcal{D}'_{ij}
$$
with $\mathcal{D}'_{ij}=\left\{ d_{ij}-Q g_j+g_i\ :\  d_{ij} \in \mathcal{D}_{ij}, 1 \le i,j \le m \right\}$, or in short,
$\mathcal{D}'_{ij}=\mathcal{D}_{ij}-Q g_j+g_i$. The corresponding tile equation becomes
$$
Q A'_j= \bigcup A'_i + \mathcal{D}'_{ij}
$$
which is satisfied by $A'_j=A_j+g_j$. So we set $\supp(T'_j)=A'_j$ and
the color of $T'_j$ to be the one of $T_j$. 
Such modification only translates 
supports of tiles and corresponding point sets, which causes 
non essential changes of the description of $\Tk$. Hereafter, we do not 
distinguish such changes of reference points and use the 
same symbols $\Lam_i$ and $T_i$. 

Since the representative points of tiles in a tiling $\Tk$ are
taken in the relatively same position for the same type of tiles, for each $i \le m$ we define a {\em reference point} $c_i \in \R^d$ for the prototile $T_i$ and can take an associated substitution Delone multi-color set $\Lb = \Lb_{\Tk} = (\Lam_i)_{i \le m}$ satisfying $\Tk = \{T_i - c_i + u \ | \ u \in \Lam_i, i \le m \}$.

The topological dynamical system $(X_{\Tk},\R^d)$ has a unique invariant Borel probability measure $\mu$ when the tiling $\Tk$ is a self-affine tiling
(see \cite{soltil, LMS2}).
The dynamical spectrum of the system $(X_{\Tk}, \R^d, \mu)$ refers to the spectrum of the unitary operator $U_x$ arising from the translational action on the space of $L^2$-functions on $X_{\Tk}$. 
We say that the tiling $\Tk$ has {\em pure discrete dynamical spectrum} if the eigenfunctions for the $\R^d$-action span a dense subspace of $L^2(X_{\Tk},\mu)$ \cite{soltil}.

\subsection{Overlap coincidence} \label{Overlap-coincidence}

Let $\Xi(\Tk) = \{y \in \R^d \ : \ T = y + S, \ \mbox{where $T, S \in \Tk$}\}$.
A triple $(u, y, v)$, with $u + T_i, v + T_j \in \Tk$
and $y \in \Xi(\Tk)$, is called an {\em overlap} if
\[ (u+A_i-y)^{\circ} \cap (v+A_j)^{\circ} \neq \emptyset, \]
where $A_i = \supp(T_i)$ and $A_j = \supp(T_j)$. An overlap $(u, y, v)$ is a {\em coincidence} if
\[ \mbox{$u-y = v$ and $u + T_i, v + T_i \in \Tk$ for some $i \le m$}.\]
Let $\mathcal{O} = (u, y, v)$ be an overlap in $\Tk$, we define
{\em $\ell$-th inflated overlap}
\begin{eqnarray*}
{\Omega}^{\ell} \mathcal{O} = \{(u', Q^{\ell} y, v') \, :
u'+T_k \in \Omega^{\ell}(u + T_i), v'+T_r \in \Omega^{\ell}(v+T_j), \ \mbox{and $(u',Q^{\ell}y,v')$ 
is an overlap} \}.
\end{eqnarray*}

\begin{defi} \label{def-overlapCoincidence}
{\em We say that a self-affine tiling $\Tk$ admits {\em overlap
coincidence} if there exists $\ell \in \Z_+$ such that for each
overlap $\mathcal{O}$ in $\Tk$, ${\Omega}^{\ell} \mathcal{O}$
contains a coincidence.}
\end{defi}

\begin{theorem} \cite{LMS2, Lee}
Let $\Tk$ be a self-affine tiling in $\R^d$ for which $\Xi(\Tk)$ is a Meyer set. Then $(X_{\Tk}, \R^d, \mu)$ has pure discrete
dynamical spectrum if and only if $\Tk$ admits overlap
coincidence.
\end{theorem}

\begin{theorem} \cite{Lee} \label{Lee-theorem}
Let $\Tk$ be a 
self-affine tiling in $\R^d$ for which $\Xi(\Tk)$ is a Meyer set and 
$\Lb = \Lb_{\Tk} = (\Lam_i)_{i \le m}$ be an associated substitution 
Delone multi-color set.
Then $\Tk$ admits overlap coincidence if and only if $\Lb$ admits algebraic coincidence.
\end{theorem}

\section{Strong coincidence and Simultaneous coincidence in $\R^d$}

A substitution $\sigma$ over $m$ letters $\{1,2,\dots,m\}$ 
is called primitive, when its substitution matrix $M_{\sigma}$
is primitive. We say that $\sigma$ is {\it irreducible} if the characteristic polynomial of $M_{\sigma}$ is irreducible, and it is {\it Pisot}, if the Perron 
Frobenius root $\beta$
of 
$M_{\sigma}$ is a Pisot number. 
Then we can construct a bi-infinite sequence which is generated by $\sigma$ from some fixed letters around the origin and define the natural {\em suspension tiling} $\Tk$ in $\R$
with an expansion factor $\beta$
arose from $\sigma$
by associating to the letters of the bi-infinite sequence the intervals whose lengths are given by a left eigenvector 
of $M_{\sigma}$ corresponding to $\beta$. In this case, $\beta$ is identified with 
the $1\times 1$ expansive matrix $Q=(\beta)$. 
Let $\Omega$ be the corresponding tile-substitution 
and let $\mathcal{A} = \{T_1, \dots, T_m\}$ be the corresponding prototiles of the intervals whose left end points are all at the origin.
Taking the left end points of the intervals of tiles in the tiling, we can get an associated substitution Delone multi-color set.
We can interpret ``strong coincidence'' of $\sigma$
into the tiling setting in the following way:

\begin{defi}
\label{1D}
{\em Let $\Tk$ be the suspension tiling in $\R$ of a substitution $\sigma$.
We say that $\Tk$ admits {\em prefix strong coincidence} if for any pair of prototiles $\{T_i, T_j\} \subset \mathcal{A}$, 
there exists $L \in \N$} such that two supertiles $\Omega^L T_i$ and $\Omega^L T_j$ have at least one common tile where $\Omega$ is the tile-substitution. 
\end{defi}

We generalize this definition of strong coincidence to $\R^d$ by using 
the reference points. 

\begin{defi} \label{admissible}
{\em Let $\Tk$ be a self-affine tiling in $\R^d$.
Let $c_i$ be the reference point of $T_i$ for $i \le m$.  
We say that the set of the reference points 
is {\it admissible} if $\cap_{i \le m} (\supp(T_i)-c_i)$ has non-empty interior. }
\end{defi}
Clearly, the left end points of the 1-dimensional suspension tiling form a set of admissible reference points.

\begin{defi}
\label{StrongTile}
{\em Let $\Tk$ be a self-affine tiling in $\R^d$ with an expansion $Q$. Let $\mathcal{A} = \{T_1, \cdots, T_m\}$ be the prototile set of $\Tk$. 
Let $c_i$ be the reference point of $T_i$ for $ i \le m$ whose set is admissible. 
Let $\Lb$ be an associated substitution Delone multi-color set
for which $\Tk = \{u_i + (T_i -c_i) \ |  \ u_i \in \Lam_i, i \le m \}$.
If for any $1 \le i, j \le m$, there is a positive integer $L$ that
\begin{eqnarray} \label{StronCoin-I}
\Omega^L(T_i - c_i) \cap 
\Omega^L(T_j - c_j) \neq \emptyset,
\end{eqnarray}
i.e. the left hand side contains at least one tile, then we say that $\Lb$ admits {\em strong coincidence}.}
\end{defi}

Set
$$
\mathcal{G} :=\bigcup_{k=0}^{\infty} Q^{-k} (\Lam_i-\Lam_i), \ \ \  \mbox{for some $i \le m$}.
$$
This set is independent of the choice of $i$. Indeed 
by the primitivity of $\Lb$, for any $i,j \le m$, there exists $n \in \N$ that 
$$
Q^n(\Lam_i-\Lam_i) \subset \Lam_j-\Lam_j \,.
$$
In plain words,
$\mathcal{G}$ is the set of eventual return vectors, i.e., 
vectors $v\in \R^d$ such that there is an
$n \in \N$ that $Q^n v$ is a return vector of $\Tk$.

\begin{remark}\label{Strong}
{\em If $\Lb$ admits strong coincidence, then notice that for any $ i,j \le m$, there is a common tile $T_k - c_k + \eta \in \Omega^L(T_i - c_i) \cap \Omega^L(T_j - c_j)$ for some $\eta \in \R^d$, $L \in \N$ and $ k \le m$.
Thus
$Q^L \Lam_i + \eta \subset \Lam_k$ and
$Q^L \Lam_j + \eta \subset \Lam_k$. 
So 
\begin{equation}
\label{StronCoin-II}
Q^L (\Lam_i \cup \Lam_j)\subset \Lam_k -\eta\,.
\end{equation}
Thus for any $i,j\le m$, 
\begin{equation} 
\label{Lami-Lamj_beReturnVector}
Q^L (\Lam_i-\Lam_j) \subset \Lam_k-\Lam_k \
\end{equation} 
and we obtain
$$
\Lam-\Lam \subset \mathcal{G}, \ \ \ \mbox{where} \ \Lam = \cup_{i \le m} \Lam_i \,.
$$ 
}
\end{remark}

\medskip

\begin{remark}
\label{Group}
{\em  
Assume that $\Lb = \Lb_{\Tk} = (\Lam_i)_{i \le m}$ admits
algebraic coincidence. There is $L\in \N$ and
$\eta\in \R^d$ that $Q^L(\Lam_i-\Lam_i) \subset \Lam_i-\eta$ for some $i$. This implies
$$
Q^L(\Lam_i-\Lam_i) -Q^L(\Lam_i-\Lam_i) \subset \Lam_i-\Lam_i.
$$
This is equivalent to the fact that the set of eventual return vectors
$\mathcal{G}$ forms an additive group, i.e., 
$\mathcal{G}=\langle \mathcal{G}\rangle_{\Z}$, where
$\langle \mathcal{G} \rangle_{\Z}$ is the additive group generated by the elements of $\mathcal{G}$.
Among self-affine tilings, it is interesting to characterize when 
$\mathcal{G}=\langle \mathcal{G}\rangle_{\Z}$ holds.
For e.g., consider the suspension tiling of Thue-Morse substitution 
$0\rightarrow 01,\ 1\rightarrow 10$ which does not admit overlap coincidence. 
Taking the left end points of the tiles in the tiling, we have
$$
\Lam_k \cap [0,\infty) = \left\{ \sum_{i=0}^{\ell} b_i 2^i\ :\ 
b_i\in \{0,1\}, \ \sum_{i=0}^{\ell} b_i\equiv k  \pmod{2} \text{ and } \ell \in \Z_{\ge 0}\right\}
$$
with $k \in \{0,1\}$. 
One can easily express each element of $\Lam_1\cap [0,\infty)$ as a
difference of $\Lam_0\cap [0,\infty)$. So we can show that
$\mathcal{G}=\Z[1/2]$, which forms a group.
This example also shows that the converse of Remark \ref{Strong} does not hold:
$\Lam-\Lam \subset \mathcal{G}$ does not imply strong coincidence.
A similar idea using beta-integer with golden mean base works in 
Example \ref{Fib2},
and we can prove that $\mathcal{G}=\Z[(1+\sqrt{5})/2]$.
We do not know yet any example of a
self-affine tiling for which $\mathcal{G}$ is not a group.}
\end{remark}

\medskip

Given a self-affine tiling $\Tk$, there are many ways to associate a substitution Delone multi-color
set. 
Indeed, Definition \ref{1D} used the left end points but one may choose
other reference points.
Arnoux-Ito \cite{AI} 
gave the strong coincidence with respect to 
right end points as well, that is, {\em 
suffix strong coincidence}\footnote{In \cite{AI}, prefix (resp. suffix) strong coincidence is called positive (resp. negative) strong coincidence.}.
There is a standard way to associate Delone multi-color set to a given self-affine tiling.
A tile map $\gamma:\Tk \rightarrow \Tk$ sends a tile $T$ to the one in $\Omega(T)$ such that $\gamma(T_1)$ and $\gamma(T_2)$ are located in the same relative position
in $\Omega(T_1)$ and $\Omega(T_2)$ whenever $T_1$ and $T_2$ have the same color.
A control point $c(T)$ of $T\in \Tk$ is defined by
$$ c(T)=\bigcap_{n=1}^{\infty} Q^{-n} (\gamma^n T).
$$
The color of the control point $c(T)$ is given by the color of the tile $(T)$.
By definition control points of the same color tiles are located in the same relative
position and the set of control points $\mathcal{C}$ is invariant under the expansion by $Q$, that is, $Q \mathcal{C} \subset \mathcal{C}$. We can choose control points
as representative points. Let $\Lam_i$ be chosen to be the control points of
tiles in $\Tk$ of color $i$. Then $\mathcal{C} = (\Lam_i)_{i \le m}$. The set of
control points $\mathcal{C}$ forms a substitution Delone multi-color set.

\begin{example} 
{\em 
\label{aba} 
The substitution $\sigma$ defined by $a\rightarrow aba, b\rightarrow bab$ satisfies
neither prefix nor suffix strong coincidence. 
However define 
the tile map $\gamma$ by choosing $T_b$ in $\sigma(T_a)$
and the left most $T_b$ in $\sigma(T_b)$, then 
$\Lb$ admits strong coincidence by the choice of the control points. 
}
\end{example}

Now we introduce a stronger notion than Definition \ref{StrongTile}:

\begin{defi}
\label{SimulTile}
{\em 
Let $\Tk$ be a self-affine tiling in $\R^d$ with expansion $Q$. Let $\mathcal{A} = \{T_1, \cdots, T_m\}$ be the prototile set of $\Tk$. 
Let $c_i$ be the reference point of $T_i$ for $i \le m$ whose set is admissible. 
Let $\Lb$ be an associated substitution Delone multi-color set
for which $\Tk = \{u_i + (T_i -c_i) \ |  \ u_i \in \Lam_i, i \le m \}$.
If there is a positive integer $L$ that
\begin{eqnarray} \label{SimulCoin-I}
\bigcap_{i=1}^m \Omega^L(T_i - c_i) \neq \emptyset,
\end{eqnarray}
i.e. the left hand side contains at least one tile, then we say that $\Lb$ admits {\em simultaneous coincidence}.
}
\end{defi}

\begin{remark}
\label{Simul}
{\em 
By the similar argument as Remark \ref{Strong},
if $\Lb$ admits simultaneous coincidence then there exists $L\in \N$,
$\eta\in \R^d$ and $k \le m$
such that 
\begin{eqnarray} \label{SimulCoin-II}
Q^L (\bigcup_{i =1}^m \Lam_i)\subset \Lam_k -\eta\,.
\end{eqnarray}
}
\end{remark}

\begin{remark}
{\em  
In $d=1$ in Arnoux-Ito's framework of \cite{AI}, they have chosen 
the left (or right) end points as control points. It is 
plausible that their definition
of geometric substitution acting on broken segments
would work by other choices of control points by shifting 
the origin to the control points located on the 1-dimensional suspension
tiling in the expanding line of $M_{\sigma}$.
In this case, a unit segment will grow 
in prefix and suffix directions at a time, which fits better with
the prefix-suffix construction by Canterini and Siegel \cite{CS}.
The definition of strong coincidence would naturally be extended in the form
to assure that some iterates of every pair of 
two unit segments share a common segment. 
}
\end{remark}

\medskip

For a Pisot substitution tiling in $\R$, one can consider a substitution Delone multi-color set $\Lb$ taking reference points from the left end points of tile intervals.
The height group of $\Lb$ is defined by Sing \cite{SingThesis} which generalizes
an idea of Dekking \cite{Dek} for constant length substitution.
We extend this definition to general substitution Delone multi-color sets in $\R^d$.

\begin{defi}
\label{Height}
{\em {\em The height group} of $\Lb = (\Lam_i)_{i \le m}$ in  $\R^d$ is the quotient group 
$$ \langle \Lam -\Lam \rangle_{\Z} \ / \ \langle \Lam_i-\Lam_i \ |\ i\le m \rangle_{\Z}, \ \ \ \mbox{where $\Lam = \cup_{i \le m}\Lam_i $}\,.$$
}
\end{defi}

Given a tiling $\Tk$ in  $\R^d$, the height group also {\bf depends} on the choice of $\Lb$. 

\begin{example} 
\label{abaH}
{\em Consider again the substitution in Example \ref{aba}. Taking the left end points of $T_a$ and $T_b$, we have 
$\Lambda=\Lam_a \cup \Lam_b = \Z$, $\langle \Lam -\Lam \rangle_{\Z}=\Z$, and $
\langle \Lam_i-\Lam_i \ |\ i \in \{a, b\} \rangle_{\Z}=2\Z$. Thus the height group is $\Z/2\Z$.
The same is true for the right end points. 
However the third choice in Example \ref{aba} gives 
the height group $(2/3)\Z/2\Z\simeq \Z/3\Z$.
We note that the height group is not trivial by any choice of
admissible reference points.
In fact, let $c_a$ and $c_b$ be reference points of the prototiles $T_a = ([0,1], a)$ and $T_b= ([0, 1], b)$ in $\R$.
To be admissible, $|c_b - c_a| < 1$. 
But if the height group were trivial, then
we have $c_b - c_a \in 2 \Z$, which shows $c_a=c_b$ and 
$\langle \Lam -\Lam \rangle_{\Z}=\Z$.
}
\end{example}


\section{From overlap coincidence to strong coincidence}

A {\em van Hove sequence} for $\R^d$ is a sequence 
$\mathcal{F}=\{F_n\}_{n \ge 1}$ of bounded measurable subsets of 
$\R^d$ satisfying
\be \label{Hove}
\lim_{n\to\infty} \Vol((\partial F_n)^{+r})/\Vol(F_n) = 0,~
\mbox{for all}~ r>0, \ \ \ \\
 \mbox{where} \ (\partial F_n)^{+r} = \{x \in \R^d : \mbox{dist} (x, F_n) \le r \} \nonumber.
\ee

For any $\mathcal{S} \subset \Tk$ and van Hove sequence $\{F_n\}_{n \ge 1}$, we define 
\[ \dens(\mathcal{S}) := \lim_{n \to \infty} \frac{\mbox{vol}(\mathcal{S} \cap F_n)}{\mbox{vol}(F_n)}\]
if the limit exists.
Here the density limit is dependent on the van Hove sequence. In the case of a self-affine tiling $\Tk$, if $\mathcal{S}$ is a set of all the translates of a patch in $\Tk$, the density exists uniformly \cite{LMS2}.

The following theorem shows that overlap coincidence of $\Tk$ implies strong coincidence of $\Lb_{\Tk}$ as long as the reference points for the tiling can be chosen to satisfy (\ref{transPoints-retVectors}). This theorem can be proved using various results for pure discrete spectrum such as \cite[Thm. 6.2]{soltil} and \cite[Prop. 1]{BSW}.

\begin{theorem} \label{OverlapC-to-SimulC}
Let $\Tk$ be a self-affine tiling in $\R^d$ 
with expansion $Q$ for which $\Xi(\Tk)$ is a Meyer set. Let $\Lb = \Lb_{\Tk} = (\Lam_i)_{i \le m}$ be an associated admissible substitution 
Delone multi-color set.
Assume that 
\be \label{transPoints-retVectors}
\Lam-\Lam \subset \langle \mathcal{G} \rangle_{\Z}, 
\ee
where $\Lam = \cup_{i \le m} \Lam_i$ and $\mathcal{G} = =\bigcup_{k=0}^{\infty} Q^{-k} (\Lam_i-\Lam_i)$ for some $i \le m$.
Then overlap coincidence of $\Tk$ implies strong coincidence of $\Lb$. Furthermore overlap coincidence of $\Tk$ implies simultaneous coincidence of $\Lb$.
\end{theorem}

\noindent {\sc Proof.}
For any $ i, j \le m$, take $v_i \in \Lam_i$ and $v_j \in \Lam_j$. Let $\alpha = v_j - v_i$. 
By the assumption, $$
\alpha \in \Lam_j-\Lam_i  \subset
\langle \mathcal{G} \rangle_{\Z}\,. $$
From Theorem \ref{Lee-theorem} and Remark \ref{Group}, $\mathcal{G} = \langle \mathcal{G} \rangle_{\Z}$.
Therefore there is $M\in \N$ such that $Q^M \alpha \in \Lam_1-\Lam_1$. Overlap coincidence implies that 
\begin{equation}
\label{dense}
{\rm dens}  \ \Tk \cap (\Tk -Q^n \alpha)  \longrightarrow 1
\end{equation}
as $n\rightarrow \infty$ \cite{soltil}.
Since $\Tk$ has uniform patch frequencies \cite{LMS2}, the density is not dependent on the choice of the van Hove sequence. So we can take a van Hove sequence $\{ \supp (\Omega^n (T_i -c_i + v_i)) \}_{n \ge 1}$. 
Then we observe that $\Tk$ admits strong coincidence, otherwise
the density is bounded away from one.

Furthermore, we claim that $\Lb$ admits simultaneous coincidence. 
Fix $k \le m$ and $v_k \in \Lam_k$. For any $i \le m$, choose $v_i \in \Lam_i$.
Let $\alpha_i = v_k - v_i$.
By the admissibility of $\Lb$, 
$\cap_{i \le m} \supp(T_i-c_i +v_k)$ has non-empty interior.
Let 
$$
P_n = \bigcap_{i \le m} \supp ( \Omega^n(T_i - c_i + v_i + \alpha_i)) \ \ \ \mbox{for $n \in \N$}.
$$
Then $\{P_n\}_{n \ge 1}$ is a van Hove sequence and
$\Omega^n (T_k-c_k +v_k)$ is a patch in $\Tk$. 
From the overlap coincidence of $\Tk$, for each $i \le m$
$$
\frac{{\rm vol}  ( \Omega^n (T_k - c_k + v_k) \cap \Omega^n ( T_i - c_i + v_i + \alpha_i) \cap P_n) }{{\rm vol} (P_n) } \ \stackrel{n \to \infty}{\longrightarrow} \ 1 \,.
$$
Thus
\begin{eqnarray*}
\lefteqn{1-\frac{{\rm vol} \left( \left( \bigcap_{i \le m, i \neq k} \Omega^n (T_k - c_k + v_k) \cap \Omega^n ( T_i - c_i + v_i + \alpha_i ) \right) \cap P_n \right)}{{\rm vol} (P_n)} } \\
& \le &
\sum_{i \le m, i \neq k} \left(
1-\frac{{\rm vol} \left( \Omega^n (T_k - c_k + v_k) \cap \Omega^n ( T_i -c_i + v_i + \alpha_i ) \cap P_n \right)}{{\rm vol} (P_n)} \right) \ \longrightarrow \ 0 \ \,. 
\end{eqnarray*}
It implies that there exists $n \in \N$ such that 
$$
\bigcap_{i \le m, i \neq k} \Omega^n (T_k - c_k + v_k) \cap \Omega^n ( T_i - c_i + v_i + \alpha_i) \cap P_n) \ \neq \ \emptyset \,.
$$
Thus 
$$
\bigcap_{i \le m, i \neq k} \Omega^n (T_k - c_k) \cap \Omega^n ( T_i - c_i) \ \neq \ \emptyset 
$$
and the claim follows.
\qed

\medskip

\begin{cor}
\label{HeightOneSimul}
Let $\Tk$ be a self-affine tiling in $\R^d$ 
for which $\Xi(\Tk)$ is a Meyer set. Let
$\Lb = \Lb_{\Tk} = (\Lam_i)_{i \le m}$ be an associated admissible substitution Delone multi-color set with a trivial height group.
Then overlap coincidence of $\Tk$ implies simultaneous coincidence as well as strong coincidence of $\Lb$.
\end{cor}

\noindent {\sc Proof.}
The trivial height group assumption implies 
$\Lam-\Lam \subset \langle \mathcal{G}
\rangle_{\Z}$.
\qed
\medskip

We provide a lemma.

\begin{lemma}
\label{LinIndep} 
Let $A$ be a $m\times m$ integer matrix whose characteristic polynomial is irreducible. 
Then the entries of an eigenvector of $A$ are linearly independent over $\Q$.
\end{lemma}

\noindent {\sc Proof.}
Let $\alpha$ be an eigenvalue of $A$ and $(v_1,v_2, \dots, v_m)\in \Q(\alpha)^m$
be a corresponding eigenvector. Denote by $\alpha_i\ (i=1,\dots, m)$
the conjugates of $\alpha$. Applying Galois conjugate map $\tau_i$
which sends $\alpha$ to $\alpha_i$, we obtain $d$ different eigenvectors 
$(\tau_i(v_1),\dots,\tau_i(v_m))$. Since $\alpha_i\ (i=1,\dots, m)$ are all distinct, 
the corresponding eigenvectors are linearly independent over $\C$. So $B: =(\tau_i(v_j))$ is invertible. 
However if $v_1,v_2, \dots, v_m$ are linearly
dependent over $\Q$, then the column vectors of $B$ are linearly dependent over $\C$ which is a contradiction.
\qed

\medskip

Now we consider an 1-dim substitution tiling and
recall the following lemma from \cite[Lemma\, 6.34]{SingThesis} and \cite[Thm. 12.1]{BK}).

\begin{lemma}\cite{SingThesis, BK}
\label{translation-group}
Let $\sigma$ be a primitive irreducible substitution and 
consider its suspension tiling of $\R$ with an expansion
factor $\beta$, the Perron Frobenius root of  the substitution matrix $M_{\sigma}$. 
Let $\Lb$ be an associated substitution Delone multi-color set whose points are taken from the left end of tiles in the suspension tiling.
Then we have
\begin{eqnarray} 
\langle \Lam_i - \Lam_i \ | \ i \le m \rangle_{\Z} = \langle (\cup_{i \le m} \Lam_i) - (\cup_{i \le m} \Lam_i) \rangle_{\Z}. 
\end{eqnarray}
\end{lemma}

\medskip

Corollary \ref{HeightOneSimul} and
Lemma \ref{translation-group} give a combinatorial result on Pisot substitution sequences: the
overlap coincidence implies
strong and simultaneous coincidence for $1$-dim irreducible Pisot substitutions.
This result agrees with the earlier results of \cite{BK, BBK}.

\begin{cor}
\label{1DSimul}
Let $\sigma$ be the irreducible Pisot substitution over letters $\{1,\dots,m\}$ 
whose natural suspension tiling satisfies overlap coincidence. 
Then there are $L\in \N$ and $M\in \N$ such that prefixes of length $M$ of
$
\sigma^L(1), \sigma^L(2), \dots, \sigma^L(m)
$
has the same number of each letter $j$ for $1 \le j \le m$,
and all ends with an identical letter.
\end{cor}

\noindent {\sc Proof.}
As in Lemma \ref{translation-group}, we choose the left end points as
the set of control points. Then it fulfills the conditions of the admissibility and the trivial height group of Theorem \ref{OverlapC-to-SimulC}.
Theorem \ref{OverlapC-to-SimulC} shows that the suspension tiling satisfies simultaneous
coincidence by taking $n$-th iterates of all prototiles. 
Since the suspension lengths of prototiles form a left eigenvector of 
$M_{\sigma}$, they are linearly independent over $\Q$ by Lemma \ref{LinIndep}. 
Counting the number of tile types in $n$-th iterate of each tile up to the tile where simultaneous coincidence occurs, we can claim the result.
\qed

\medskip

\begin{example}
\label{Fib2}
{\em We give an example of a substitution tiling for which an associated admissible substitution Delone multi-color set has a trivial height group, but it admits neither overlap coincidence nor strong coincidence (simultaneous coincidence).
Consider the substitution $\sigma$ over four letters $\{a,b,A,B\}$ defined by 
$$a\rightarrow aB, b\rightarrow a,
A\rightarrow Ab, B\rightarrow A.$$ This is $\Z/2\Z$-extension of Fibonacci
substitution. The suspension tiling dynamics does not satisfy 
overlap coincidence as it is computed in \cite{AL}. 
The suspension lengths of capital letter and non capital letter are the same. 
Introduce a letter to the letter involution $\tau$ which interchanges capital 
letter to non capital letter and vice versa. 
Then $\tau$ and $\sigma$ commutes
from $\tau(\sigma^n(a))=\sigma^n(A)$ for all $n \in \N$, we see that $\sigma$ does not 
satisfy prefix (nor suffix) strong coincidence. 
However taking the left end points as reference points, the 
height group is trivial which is confirmed by collecting all return vectors.
}
\end{example}

\begin{example}
\label{Rauzy2}
{\em Consider the substitution $\sigma$ over six letters $\{a,b,c,A,B,C\}$ defined by 
$$a\rightarrow aB, b\rightarrow aC, c\rightarrow a,
A\rightarrow Ab, B\rightarrow Ac, C\rightarrow A.$$ This is $\Z/2\Z$-extension 
of the Rauzy substitution. 
The suspension tiling dynamics of $\sigma$ satisfies overlap coincidence
which is checked by the algorithm in \cite{AL}. 
By the same reasoning as Example \ref{Fib2}, this substitution does 
not satisfy prefix nor suffix strong coincidence.
Taking the left end points and collecting all return vectors, 
we can compute that the height group is equal to $\Z/2\Z$. 

On the other hand,
taking the tile map $\gamma$ which sends $$
a\rightarrow B, b \rightarrow C, c\rightarrow a,
A\rightarrow A, B \rightarrow A, c\rightarrow A,
$$
the set of control points is admissible and satisfies simultaneous coincidence. 
One can also check that the height group associate to this control points
is trivial. So this gives an example of Corollary \ref{HeightOneSimul}
for a reducible substitution. 
}
\end{example}


\begin{remark}
{\em Nakaishi \cite{Nakaishi} claimed that unimodular Pisot substitutions with prefix 
strong coincidence 
generate pure discrete spectrum of the dynamical system. It would mean that strong coincidence implies overlap coincidence in the associated substitution tilings in $\R$. 
So together with Corollary \ref{1DSimul} the equivalence between overlap coincidence and strong coincidence could be established. 
We hope directly to observe the result of \cite{Nakaishi}, 
the implication from strong coincidence to overlap coincidence in the 
framework of this paper. 
}
\end{remark}

\section*{Acknowledgment}

\noindent
The authors would like to thank the referees for the valuable comments and references.
This research was supported by Basic Research Program through the
National Research Foundation of Korea(NRF) funded by the Ministry
of Education, Science and Technology(2010-0011150) and the
Japanese Society for the Promotion of Science (JSPS), Grant in aid
21540010. The second author is grateful for the support of KIAS for this research.

\end{document}